\documentclass{amsart}
\usepackage{graphicx}
\usepackage{amscd}
\usepackage{amsmath}
\usepackage{amsfonts}
\usepackage{amssymb}
\newtheorem{theorem}{Theorem}
\theoremstyle{plain}
\newtheorem{acknowledgement}{Acknowledgement}

\newtheorem{lemma}{Lemma}

\newtheorem{proposition}{Proposition}

\numberwithin{equation}{section}

\begin{document}
\title[Coulhon type inequalities]{A note on Coulhon type inequalities}
\author{Joaquim Mart\'{i}n$^{\ast}$}
\address{Department of Mathematics\\
Universitat Aut\`{o}noma de Barcelona}
\email{jmartin@mat.uab.cat}
\author{Mario Milman**}
\address{Department of Mathematics\\
Florida Atlantic University\\
Boca Raton, Fl. 33431}
\email{extrapol@bellsouth.net}
\urladdr{http://www.math.fau.edu/milman}
\thanks{$^{\ast}$Partially supported in part by Grants MTM2010-14946, MTM-2010-16232.}
\thanks{**This work was partially supported by a grant from the Simons Foundation
(\#207929 to Mario Milman).}
\thanks{This paper is in final form and no version of it will be submitted for
publication elsewhere.}
\subjclass{2000 Mathematics Subject Classification Primary: 46E30, 26D10.}
\keywords{Sobolev inequalities, modulus of continuity, symmetrization, isoperimetric
inequalities, interpolation.}

\begin{abstract}
T. Coulhon introduced an interesting reformulation of the usual Sobolev
inequalities. We characterize Coulhon type inequalities in terms of
rearrangement inequalities.
\end{abstract}\maketitle

\section{Introduction}

Let $(X,d,\mu)$ be a connected Borel metric measure space. The perimeter or
Minkowski content of a Borel set $A\subset X,$ is defined by
\[
\mu^{+}(A)=\lim\inf_{h\rightarrow0}\frac{\mu\left(  A_{h}\right)  -\mu\left(
A\right)  }{h},
\]
where $A_{h}=\left\{  x\in\Omega:d(x,A)<h\right\}  ,$ and the isoperimetric
profile $I=I_{(\Omega,d,\mu)}$ is defined by
\[
I_{(\Omega,d,\mu)}(t)=\inf_{A}\{\mu^{+}(A):\mu(A)=t\}.
\]
We assume throughout that $(X,d,\mu)$ is such that $I_{(\Omega,d,\mu)}$ is
concave, continuous with $I(0)=0.$ Moreover, we also assume that $(X,d,\mu)$
is such that for each $c\in R,$ and each $f\in Lip_{0}(X),|\nabla
f(x)|=0,a.e.$ in the set $\{x:f(x)=c\}.$ Under these conditions\footnote{In
\cite{mamiadv} the result is shown for metric probability spaces such that
$I(t)$ is symmetric about $1/2,$ in which case we can replace $Lip_{0}(X)$ by
$Lip(X)$ in the statement. With minor modifications one can also show its
validity for infinite measure spaces (cf. also \cite{mmp}, \cite{mamicon}).}
we showed in \cite{mamiadv} that the Gagliardo-Nirenberg-Ledoux inequality
\begin{equation}
\int_{0}^{\infty}I(\mu_{f}(t))dt\leq\left\|  \left|  \nabla f\right|
\right\|  _{L{^{1}(X)}},\text{ for all }f\in Lip_{0}(X) \label{gn}%
\end{equation}
is equivalent to
\begin{equation}
f^{\ast\ast}(t)-f^{\ast}(t)\leq\frac{t}{I(t)}\left|  \nabla f\right|
^{\ast\ast}(t), \label{c0}%
\end{equation}
where $Lip_{0}(X)$ are the functions in $Lip(X)$ of compact support,
\[
\left|  \nabla f(x)\right|  =\lim\sup_{y\rightarrow x}\frac{\left|
f(x)-f(y)\right|  }{d(x,y)},
\]
$\mu_{f}(t)=\mu\{\left|  f\right|  >t\},$ $f^{\ast}$ is the non increasing
rearrangement\footnote{For background we refer to \cite{bs} (on
rearrangements), and \cite{leoni}, \cite{Maz} (on Sobolev spaces) .} of $f$
with respect to the measure $\mu$ and $f^{\ast\ast}(t)=\frac{1}{t}\int_{0}%
^{t}f^{\ast}(s)ds.$

Conversely, if an inequality of the form (\ref{c0}) holds for some continuous
concave function $I_{1}(t)$, it was shown in \cite{mamiadv} that $I_{1}(t)$
satisfies the isoperimetric inequality $I_{1}(\mu(A))\leq\mu^{+}(A) $ for any
Borel set\footnote{Therefore, $I_{1}(t)\leq\inf\{\mu^{+}(A):\mu(A)=t\}=I(t),$
and consequently $\frac{t}{I(t)}\leq\frac{t}{I_{1}(t)}.$} $A\subset\subset X$.
In particular, for $\mathbb{R}^{n}$ it is well known that (cf. \cite[Chapter
1]{Maz}) $I(t)=c_{n}t^{1-1/n},$ and therefore (\ref{c0}) becomes (cf.
\cite{bmr} and the references therein)
\begin{equation}
f^{\ast\ast}(t)-f^{\ast}(t)\leq c_{n}^{-1}t^{1/n}\left|  \nabla f\right|
^{\ast\ast}(t). \label{berta}%
\end{equation}
It follows that (\ref{gn}) gives
\[
c_{n}\int_{0}^{\infty}\mu_{f}(t)^{1-1/n}dt=c_{n}\frac{1}{n^{\prime}}\int
_{0}^{\infty}t^{1/n^{\prime}}f^{\ast}(t)\frac{dt}{t}\leq\left\|  \left|
\nabla f\right|  \right\|  _{L{^{1}(\mathbb{R}^{n})}}%
\]
i.e.
\[
\left\|  f\right\|  _{L^{^{\frac{n}{n-1}},1}(\mathbb{R}^{n})}\leq c\left\|
\left|  \nabla f\right|  \right\|  _{L{^{1}(\mathbb{R}^{n})}},\text{ for all
}f\in Lip_{0}(\mathbb{R}^{n}).
\]
In other words, (\ref{gn}) represents a generalization of the sharp form of
the Euclidean Gagliardo-Nirenberg inequality that uses Lorentz spaces (cf.
\cite{poor} and \cite{mmp} (for Euclidean spaces), \cite{le} (Gaussian
spaces), and \cite{boho}, \cite{mamiadv} (for metric spaces); for the
corresponding rearrangement inequalities we refer to \cite{bmr}, \cite{mmjfa},
\cite{mamiadv}, as well as the references therein).

The corresponding Sobolev inequalities when $\left|  \nabla f\right|  \in
L^{p},$ $p>1$ are also known to self improve (cf. \cite{Maz}, \cite{bakr},
\cite{mamicon}, and the references therein) but an analogous rearrangement
inequality characterization in this case has remained an open problem. On the
other hand, Coulhon (cf. \cite{cou1}, \cite{cou2}, \cite{cou}) and
Bakry-Coulhon-Ledoux \cite{bakr} introduced and studied a different scale of
Sobolev inequalities. For $p\in\lbrack1,\infty],$ and $\phi$ an increasing
function on the positive half line, these authors studied the validity of
inequalities of the form
\[
(S_{\phi}^{p})\;\;\;\;\left\|  f\right\|  _{p}\leq\phi(\left\|  f\right\|
_{0})\left\|  \left|  \nabla f\right|  \right\|  _{p},\text{ }f\in
Lip_{0}(X),
\]
where
\[
\left\|  f\right\|  _{0}=\mu\{support(f)\},\text{ }\left\|  f\right\|
_{p}=\left\|  f\right\|  _{L^{p}(X)}.
\]

In particular, it was shown by Coulhon et al. that the $(S_{\phi}^{p})$
inequalities encapsulate the classical Sobolev inequalities, as well as the
Faber-Krahn inequalities. For $p=1,$ $(S_{\phi}^{1})$ is equivalent to the
isoperimetric inequality in the sense that\footnote{See Section \ref{sec1}%
\ below.}
\[
\frac{t}{I(t)}\leq\phi(t).
\]
Moreover, for $p=\infty,$ the $(S_{\phi}^{\infty})$ conditions are explicitly
connected with volume growth. For a detailed discussion of the different
geometric interpretations for different $p^{\prime}s$ we refer to \cite{cou1},
\cite{gri}, \cite{Maz}, and the references quoted therein.

It follows from this discussion that, for a suitable class of metric measure
spaces, the $(S_{\phi}^{1})$ condition can be characterized by means of the
symmetrization inequality (\ref{c0}):
\[
(S_{\phi}^{1})\text{ holds }\Leftrightarrow\text{(\ref{c0}) holds.}%
\]
The purpose of this paper is to provide an analogous rearrangement
characterization of the $(S_{\phi}^{p})$ conditions, \ for $1\leq p<\infty.$
Our main result extends (\ref{c0}) as follows

\begin{theorem}
\label{teo}Let $(X,d,\mu)$ be a connected Borel metric measure space as
described above, and let $p\in\lbrack1,\infty).$ The following statements are equivalent

\begin{enumerate}
\item $(S_{\phi}^{p})$ holds, i.e.
\begin{equation}
\left\|  f\right\|  _{p}\leq\phi(\left\|  f\right\|  _{0})\left\|  \left|
\nabla f\right|  \right\|  _{p},\text{ for all }f\in Lip_{0}(X).
\label{desigualdad01}%
\end{equation}

\item  Let $k\in\mathbb{N}$ be such that $k<p\leq k+1,$ then for all $f\in
Lip_{0}(X)$
\begin{equation}
\left(  \frac{f_{\left(  p\right)  }^{\ast\ast}(t)}{\phi_{(p)}(t)}\right)
^{1/p}-\left(  \frac{f_{(p)}^{\ast}(t)}{\phi_{(p)}(t)}\right)  ^{1/p}%
\leq2^{\frac{k+1}{p}-1}\left(  \left|  \nabla f\right|  _{(p)}^{\ast\ast
}(t)\right)  ^{1/p}, \label{norma}%
\end{equation}
where
\[
f_{(p)}^{\ast}(t)=\left(  f^{\ast}(t)\right)  ^{p},\text{ }f_{(p)}^{\ast\ast
}(t)=\frac{1}{t}\int_{0}^{t}f_{(p)}^{\ast}(s)ds,\text{ }\phi_{(p)}(t)=\left(
\phi(t)\right)  ^{p}.
\]

\item  Let $k\in\mathbb{N}$ be such that $k<p\leq k+1,$ then for all $f\in
Lip_{0}(X),$ $f_{(p)}^{\ast}$ is absolutely continuous (cf. \cite{leoni}) and
\begin{equation}
-\frac{\partial}{\partial t}\left(  f_{(p)}^{\ast\ast}(t)\right)
^{1/p}=-\frac{\partial}{\partial t}\left(  \frac{1}{t}\int_{0}^{t}%
f_{(p)}^{\ast}(s)ds\right)  ^{1/p}\leq2^{\frac{k+1}{p}}\frac{\phi(t)}%
{t}\left(  \left|  \nabla f\right|  _{(p)}^{\ast\ast}(t)\right)  ^{\frac{1}%
{p}}. \label{norma01}%
\end{equation}
\end{enumerate}
\end{theorem}

Note that for $p=1$ the inequality (\ref{norma}) of Theorem \ref{teo}
coincides with (\ref{c0}). This new characterization for $p\geq1$ is
independent of \cite{mamiadv}, and, in fact, it provides a new approach to
(\ref{c0}) as well. On the other hand, as it is well known (cf. \cite{cou1}),
the $(S_{\phi}^{p})$ conditions get progressively weaker as $p$ increases.
Indeed, below we will also show that (\ref{c0}) implies (\ref{norma}) via an
extended form of the chain rule, that is valid for metric spaces.

The note is organized as follows. In section \ref{sec1} we give a somewhat
more detailed discussion of the $(S_{\phi}^{p})$ conditions and, in
particular, we develop a connection with \cite{mamiadv}. In section
\ref{secproof} we provide a proof of Theorem \ref{teo} and, finally, in
section \ref{secrem}, we discuss, rather briefly, connections with Nash type
inequalities, Sobolev and Faber-Krahn inequalities and
interpolation/extrapolation theory.

As usual, the symbol $f\simeq g$ will indicate the existence of a universal
constant $c>0$ (independent of all parameters involved) so that $(1/c)f\leq
g\leq c\,f$, while the symbol $f\preceq g$ means that $f\leq c\,g$. \ \ \ \ 

\section{The $(S_{\phi}^{p})$ conditions\label{sec1}}

From now on $(X,d,\mu)$ will be a connected metric measure space with a
continuous isoperimetric profile $I$ such that $\frac{t}{I(t)}$ increases and
such that $I(0)=0.$ Moreover, we also assume that $(X,d,\mu)$ is such that for
each $c\in\mathbb{R},$ and each $f\in Lip_{0}(X),|\nabla f(x)|=0,a.e.$ in the
set $\{x:f(x)=c\}.$ The isoperimetric profile $I=I_{(\Omega,d,\mu)}$ is
defined by
\[
I_{(\Omega,d,\mu)}(t)=\inf_{A}\{\mu^{+}(A):\mu(A)=t\},
\]
where $\mu^{+}(A)$ is the perimeter or Minkowski content of the Borel set
$A\subset X,$ defined by
\[
\mu^{+}(A)=\lim\inf_{h\rightarrow0}\frac{\mu\left(  A_{h}\right)  -\mu\left(
A\right)  }{h},
\]
where $A_{h}=\left\{  x\in\Omega:d(x,A)<h\right\}  .$

\subsection{The $(S_{\phi}^{1})$ condition}

From \cite{mamiadv} (cf. also \cite{mamicon}) we know that
\begin{equation}
f^{\ast\ast}(t)-f^{\ast}(t)\leq\frac{t}{I(t)}\left|  \nabla f\right|
^{\ast\ast}(t),\text{ }f\in Lip_{0}(X), \label{c2.1}%
\end{equation}
is equivalent to the isoperimetric inequality. If we combine these results
with the characterization of $(S_{\phi}^{1})$ given in \cite{cou} we can see
the equivalence between (\ref{c2.1}) and the $(S_{\phi}^{1})$ condition. To
understand the discussion of the next section it is instructive to provide an
elementary direct approach. So we shall now show that (\ref{c2.1}) implies
$(S_{\phi}^{1})$ with $\phi(t)=t/I(t),$ and that this choice is in some sense
the best possible $(S_{\phi}^{1})$ condition.

Suppose that (\ref{c2.1}) holds. Multiplying both sides of (\ref{c2.1}) by
$t>0$ we obtain
\[
t\left(  f^{\ast\ast}(t)-f^{\ast}(t)\right)  \leq\frac{t}{I(t)}\int_{0}%
^{t}\left|  \nabla f\right|  ^{\ast}(s)ds.
\]
Since formally $f^{\ast}(t)=\mu_{f}^{-1}(t),$ drawing a diagram it is easy to
convince oneself that
\begin{align*}
t\left(  f^{\ast\ast}(t)-f^{\ast}(t)\right)   &  =\int_{0}^{t}f^{\ast
}(s)ds-tf^{\ast}(t)\\
&  =\int_{f^{\ast}(t)}^{\infty}\mu_{f}(s)ds.
\end{align*}
Consequently, if we let $t=\left\|  f\right\|  _{0},$ we see that $f^{\ast
}(\left\|  f\right\|  _{0})=0,$ $\int_{f^{\ast}(\left\|  f\right\|  _{0}%
)}^{\infty}\mu_{f}(s)ds=\left\|  f\right\|  _{1},$ and $\int_{0}^{\left\|
f\right\|  _{0}}\left|  \nabla f\right|  ^{\ast}(s)ds=\left\|  \left|  \nabla
f\right|  \right\|  _{1}.$ Thus,
\[
\left\|  f\right\|  _{1}\leq\frac{\left\|  f\right\|  _{0}}{I(\left\|
f\right\|  _{0})}\left\|  \left|  \nabla f\right|  \right\|  _{1}.
\]
In other words, the $(S_{\phi}^{1})$ condition holds with $\phi(t)=\frac
{t}{I(t)},$ and consequently the $(S_{\tilde{\phi}}^{1})$ condition holds for
any $\tilde{\phi}(t)\geq\frac{t}{I(t)}.$ On the other hand, consider an
$(S_{\tilde{\phi}}^{1})$ condition for a continuous, increasing but arbitrary
function $\tilde{\phi}$. Let $A$ be a Borel set, $A\subset\subset X,$ with
$\mu(A)=t.$ Formally inserting $f=\chi_{A}$ in the corresponding
$(S_{\tilde{\phi}}^{1})$ inequality (this is done rigorously by
approximation), yields
\[
\left\|  \chi_{A}\right\|  _{1}=t=\mu(A)\leq\tilde{\phi}(t)\mu^{+}(A).
\]
Consequently,
\begin{align*}
\frac{t}{\tilde{\phi}(t)}  &  \leq\inf\{\mu^{+}(B):\mu(B)=t\}\\
&  =I(t),
\end{align*}
and therefore
\[
\frac{t}{I(t)}\leq\tilde{\phi}(t)\text{. }%
\]

\subsection{$(S_{\phi}^{1})\Rightarrow(S_{\phi}^{p}),$ $p>1$}

In the Euclidean space $\mathbb{R}^{n}$, $I(t)=d_{n}t^{1-1/n},$ $\phi(t)\simeq
t^{1/n}$ and the best possible $(S_{\phi}^{1})$ inequality can be written as
\[
\left\|  f\right\|  _{1}\leq c_{n}\left\|  f\right\|  _{0}^{1/n}\left\|
\left|  \nabla f\right|  \right\|  _{1}.
\]

As was shown in \cite{cou1} the corresponding inequalities for $p>1$ then
follow by the (classical) chain rule, the fact that $\left\|  \left|
f\right|  ^{p}\right\|  _{0}=\left\|  \left|  f\right|  \right\|
_{0}=\left\|  f\right\|  _{0},$ and H\"{o}lder's inequality. In detail,
\begin{align*}
\left\|  f\right\|  _{p}^{p}  &  =\left\|  \left|  f\right|  ^{p}\right\|
_{1}\\
&  \leq c_{n}p\left\|  f\right\|  _{0}^{1/n}\left\|  \left|  f\right|
^{p-1}\left|  \nabla\left|  f\right|  \right|  \right\|  _{1}\\
&  \leq c_{n}p\left\|  f\right\|  _{0}^{1/n}\left\|  f\right\|  _{p}%
^{p-1}\left\|  \left|  \nabla\left|  f\right|  \right|  \right\|  _{p}.
\end{align*}
Consequently,
\[
\left\|  f\right\|  _{p}\leq c_{n}p\left\|  f\right\|  _{0}^{1/n}\left\|
\left|  \nabla f\right|  \right\|  _{p},
\]
and therefore, modulo constants, we have that $(S_{\phi}^{1})\Rightarrow
(S_{\phi}^{p}),$ for $p>1.$ More generally, this argument, taken from
\cite{cou1}, shows that the $(S_{\phi}^{p})$ conditions become weaker as $p$
increases. In the general setting of metric spaces, the classical chain rule
needs to be replaced by an inequality\footnote{The underlying elementary
inequality is
\[
\left|  a^{r}-b^{r}\right|  \leq r\left|  a^{r-1}+b^{r-1}\right|  \left|
a-b\right|  .
\]
}: for $r>1,$%
\begin{equation}
\left|  \nabla f^{r}(x)\right|  \leq2r\left|  f^{r-1}(x)\right|  \left|
\nabla f(x)\right|  . \label{chain}%
\end{equation}

Next, we use the generalized chain rule to explain the origin of the awkward
looking condition (\ref{norma}). Informally, we shall now show
that\footnote{With slightly more labor the same method will similarly show
that, more generally, $(S_{\phi}^{p})\Rightarrow(S_{\phi}^{q}),$ for $q>p.$}
$(S_{\phi}^{1})\Rightarrow(S_{\phi}^{p})$ at the level of rearrangements, i.e.
(\ref{c0})$\Rightarrow$(\ref{norma}).

Assume the validity of $(S_{\phi}^{1}).$ Let $f\in Lip_{0}(X);$ we may assume
without loss that $f$ is positive. Apply the $(S_{\phi}^{1})$ inequality to
$f_{(p)}=f^{p},$ where $p>1$ is fixed. Then, by the chain rule (\ref{chain})
\begin{align*}
f_{(p)}^{\ast\ast}(t)-f_{(p)}^{\ast}(t)  &  \preceq\phi(t)\left|  \nabla
f\right|  _{(p)}^{\ast\ast}(t)\\
&  \preceq\phi(t)(f^{p-1}\left|  \nabla f\right|  )^{\ast\ast}(t).
\end{align*}
By a result due to O'Neil (cf. \cite[page 88, Exercise 10]{bs}) and
H\"{o}lder's inequality
\begin{align*}
(f^{p-1}\left|  \nabla f\right|  )^{\ast\ast}(t)  &  \leq\frac{1}{t}\int
_{0}^{t}(f^{\ast}(s))^{p-1}\left|  \nabla f\right|  ^{\ast}(s)ds\\
&  \leq\frac{1}{t}\left(  \int_{0}^{t}f_{(p)}^{\ast}(s)ds\right)
^{1/p^{\prime}}\left(  \int_{0}^{t}\left|  \nabla f\right|  _{(p)}^{\ast
}(s)ds\right)  ^{1/p}\\
&  =\left(  f_{(p)}^{\ast\ast}(t)\right)  ^{1-1/p}\left(  \left|  \nabla
f\right|  _{(p)}^{\ast\ast}(t)\right)  ^{1/p}.
\end{align*}
Combining inequalities we obtain,
\[
f_{(p)}^{\ast\ast}(t)-f_{(p)}^{\ast}(t)\preceq\phi(t)\left(  f_{(p)}^{\ast
\ast}(t)\right)  ^{1-1/p}\left(  \left|  \nabla f\right|  _{(p)}^{\ast\ast
}(t)\right)  ^{1/p}.
\]
Hence,
\[
\left(  f_{(p)}^{\ast\ast}(t)\right)  ^{1/p}-\frac{f_{(p)}^{\ast}(t)}{\left(
\left|  f\right|  _{(p)}^{\ast\ast}(t)\right)  ^{1/p^{\prime}}}\preceq
\phi(t)\left(  \left|  \nabla f\right|  _{(p)}^{\ast\ast}(t)\right)  ^{1/p}.
\]
But, since
\[
\left(  \left|  f\right|  _{(p)}^{\ast\ast}(t)\right)  ^{1/p^{\prime}}%
\geq\left(  \left|  f\right|  _{(p)}^{\ast}(t)\right)  ^{1/p^{\prime}}=\left(
\left|  f\right|  _{(p)}^{\ast}(t)\right)  ^{1-1/p},
\]
we have
\[
\left(  f_{(p)}^{\ast}(t)\right)  ^{1/p}\geq\frac{f_{(p)}^{\ast}(t)}{\left(
\left|  f\right|  _{(p)}^{\ast\ast}(t)\right)  ^{1/p^{\prime}}},
\]
and we conclude that
\begin{align*}
\left(  f_{(p)}^{\ast\ast}(t)\right)  ^{1/p}-\left(  f_{(p)}^{\ast}(t)\right)
^{1/p}  &  \preceq\left(  f_{(p)}^{\ast\ast}(t)\right)  ^{1/p}-\frac
{f_{(p)}^{\ast}(t)}{\left(  \left|  f\right|  _{(p)}^{\ast\ast}(t)\right)
^{1/p^{\prime}}}\\
&  \preceq\phi(t)\left(  \left|  \nabla f\right|  _{(p)}^{\ast\ast}(t)\right)
^{1/p}.
\end{align*}
Therefore,
\[
\left(  \frac{f_{(p)}^{\ast\ast}(t)}{\phi_{(p)}(t)}\right)  ^{1/p}-\left(
\frac{f_{(p)}^{\ast}(t)}{\phi_{(p)}(t)}\right)  ^{1/p}\preceq\left(  \left|
\nabla f\right|  _{(p)}^{\ast\ast}(t)\right)  ^{1/p},
\]
and (\ref{norma}) holds.

\section{Proof of Theorem \ref{teo}\label{secproof}}

Before going through the proof let us make a few useful remarks. Let $\left[
x\right]  _{+}=\max(x,0),$ and let $f\geq0,$ then, for all $\lambda>0$, we
have
\begin{align}
\int_{\{f>\lambda\}}\left(  f(s)-\lambda\right)  d\mu(s)  &  =\int
[f(s)-\lambda]_{+}\,d\mu(s)=\int_{0}^{\infty}[f^{\ast}(s)-\lambda
]^{+}\,ds\label{cata}\\
&  =\int_{0}^{\infty}\mu_{\lbrack f^{\ast}-\lambda]_{+}}(s)\,ds=\int_{\lambda
}^{\infty}\mu_{f^{\ast}}(s)\,ds=\int_{\lambda}^{\Vert f\Vert_{\infty}}\mu
_{f}(s)\,ds.\nonumber
\end{align}
Thus, inserting $\lambda=f^{\ast}(t)$ in (\ref{cata}), and taking into account
that $f^{\ast}$ is decreasing, we obtain
\begin{align*}
t(f^{\ast\ast}(t)-f^{\ast}(t))  &  =\int_{0}^{t}(f^{\ast}(x)-f^{\ast
}(t))\,dx=\int_{0}^{\infty}[f^{\ast}(x)-f^{\ast}(t)]_{+}\,dx\\
&  =\int_{\{f>f^{\ast}(t)\}}\left[  f(s)-f^{\ast}(t)\right]  _{+}d\mu(s).
\end{align*}
In order to deal with $L^{p}$ norms, $p>1,$ we need to extended the formulae
above. This will be achieved through the following variant of the binomial
formula, whose proof will be provided at the end of this section.

\begin{lemma}
\label{des}Let $p>1,$ and let $k\in\mathbb{N}$ be such that $k<p\leq k+1$.
Then, for\ $a\geq b\geq0,$
\begin{equation}
(a-b)^{p}\geq a^{p}-b^{p}-\sum_{j=1}^{k}\left(
\begin{array}
[c]{c}%
p\\
j
\end{array}
\right)  b^{p-j}(a-b)^{j}, \label{des1}%
\end{equation}
and
\begin{equation}
a^{p}+b^{p}+\sum_{j=1}^{k}\left(
\begin{array}
[c]{c}%
p\\
j
\end{array}
\right)  b^{p-j}(a-b)^{j}\leq(c(p)a+b)^{p}, \label{des2}%
\end{equation}
where $c(p)=2^{\frac{k+1}{p}-1}.$
\end{lemma}

We are now ready to give the proof of Theorem \ref{teo}.

\begin{proof}
$1\rightarrow2.$ Suppose that $(S_{\phi}^{p})$ holds. We may assume without
loss that $f$ is positive. Let $t>0;$ we will apply (\ref{desigualdad01}) to
$\left[  f-f^{\ast}(t)\right]  _{+}$. Observe that
\[
\left\|  \left[  f-f^{\ast}(t)\right]  _{+}\right\|  _{0}=\mu\{f>f^{\ast
}(t)\}\leq t,
\]
and, moreover, since $\int_{\{f=f^{\ast}(t)\}}\left|  \nabla\left[
f(x)-f^{\ast}(t)\right]  \right|  dx=0,$
\[
\left\|  \nabla\left[  f-f^{\ast}(t)\right]  _{+}\right\|  _{L^{p}}^{p}%
=\int_{\{f>f^{\ast}(t)\}}\left(  \left|  \nabla f\right|  ^{\ast}(s)\right)
^{p}ds.
\]
Therefore,
\begin{align}
\left\|  \left[  f-f^{\ast}(t)\right]  _{+}\right\|  _{p}^{p}  &  \leq\left\{
\phi(\left\|  \left[  f-f^{\ast}(t)\right]  _{+}\right\|  _{0})\right\}
^{p}\left\|  \nabla\left[  f-f^{\ast}(t)\right]  _{+}\right\|  _{L^{p}}%
^{p}\nonumber\\
&  \leq\phi(t)^{p}\int_{\{f>f^{\ast}(t)\}}\left(  \left|  \nabla f\right|
^{\ast}(s)\right)  ^{p}ds\nonumber\\
&  \leq t\phi(t)^{p}\left(  \frac{1}{t}\int_{0}^{t}\left(  \left|  \nabla
f\right|  ^{\ast}(s)\right)  ^{p}ds\right) \nonumber\\
&  =t\phi(t)^{p}\left|  \nabla f\right|  _{(p)}^{\ast\ast}(t). \label{Sk}%
\end{align}
Now,
\begin{align}
\left\|  \left[  f-f^{\ast}(t)\right]  _{+}\right\|  _{p}^{p}  &
=\int_{\{f>f^{\ast}(t)\}}\left(  f(s)-f^{\ast}(t)\right)  ^{p}d\mu
(s)\nonumber\\
&  \geq\int_{\{f>f^{\ast}(t)\}}\left(  f^{p}(s)-f^{\ast}(t)^{p}\right)
d\mu(s)\nonumber\\
&  -\sum_{j=1}^{k}\left(
\begin{array}
[c]{c}%
p\\
j
\end{array}
\right)  f^{\ast}(t)^{p-j}\int_{\{f>f^{\ast}(t)\}}(f(s)-f^{\ast}(t))^{j}%
d\mu(s)\text{ \ (by (\ref{des1}))}\nonumber\\
&  =\int_{\{f_{(p)}>f_{(p)}^{\ast}(t)\}}\left(  f_{(p)}(s)-f_{(p)}^{\ast
}(t)\right)  d\mu(s)\nonumber\\
&  -\sum_{j=1}^{k}\left(
\begin{array}
[c]{c}%
p\\
j
\end{array}
\right)  f^{\ast}(t)^{p-j}\int_{\{f>f^{\ast}(t)\}}(f(s)-f^{\ast}(t))^{j}%
d\mu(s)\nonumber\\
&  =t\left(  f_{(p)}^{\ast\ast}(t)-f_{(p)}^{\ast}(t)\right)  -\sum_{j=1}%
^{k}\left(
\begin{array}
[c]{c}%
p\\
j
\end{array}
\right)  f^{\ast}(t)^{p-j}\int_{\{f>f^{\ast}(t)\}}(f(s)-f^{\ast}(t))^{j}%
d\mu(s). \label{paso1}%
\end{align}
We estimate each of the integrals in the sum using H\"{o}lder's inequality as
follows,
\begin{align}
\int_{\{f>f^{\ast}(t)\}}(f(s)-f^{\ast}(t))^{j}d\mu(s)  &  \leq\left(
\int_{\{f>f^{\ast}(t)\}}(f(s)-f^{\ast}(t))^{p}d\mu(s)\right)  ^{\frac{j}{p}%
}\left(  \int_{\{f>f^{\ast}(t)\}}d\mu(s)\right)  ^{\frac{p-j}{p}}\nonumber\\
&  =\left(  \int_{\{f>f^{\ast}(t)\}}(f(s)-f^{\ast}(t))^{p}d\mu(s)\right)
^{\frac{j}{p}}\left(  \mu_{f}(f^{\ast}(t))\right)  ^{\frac{p-j}{p}}\nonumber\\
&  \leq\left(  \int_{\{f>f^{\ast}(t)\}}(f(s)-f^{\ast}(t))^{p}d\mu(s)\right)
^{\frac{j}{p}}t^{\frac{p-j}{p}}\nonumber\\
&  =\left\|  \left[  f-f^{\ast}(t)\right]  _{+}\right\|  _{p}^{j}t^{\frac
{p-j}{p}}\nonumber\\
&  \leq\phi(t)^{j}\left(  \left|  \nabla f\right|  _{(p)}^{\ast\ast
}(t)\right)  ^{\frac{j}{p}}t^{\frac{j}{p}}t^{\frac{p-j}{p}}\text{ \ \ (by
(\ref{Sk}))}\nonumber\\
&  =t\phi(t)^{j}\left(  \left|  \nabla f\right|  _{(p)}^{\ast\ast}(t)\right)
^{\frac{j}{p}}. \label{paso2}%
\end{align}
Combining (\ref{paso1}) and (\ref{paso2}) we get
\begin{align*}
\left\|  \left[  f-f^{\ast}(t)\right]  _{+}\right\|  _{p}^{p}  &  =t\left(
f_{(p)}^{\ast\ast}(t)-f_{(p)}^{\ast}(t)\right)  -\sum_{j=1}^{p-1}\left(
\begin{array}
[c]{c}%
p\\
j
\end{array}
\right)  f^{\ast}(t)^{p-j}\int_{\{f>f^{\ast}(t)\}}(f^{\ast}(t)-f(s))^{j}%
d\mu(s)\\
&  \geq t\left(  f_{(p)}^{\ast\ast}(t)-f_{(p)}^{\ast}(t)\right)  -t\left(
\sum_{j=1}^{p-1}\left(
\begin{array}
[c]{c}%
p\\
j
\end{array}
\right)  f^{\ast}(t)^{p-j}\phi(t)^{j}\left(  \left|  \nabla f\right|
_{(p)}^{\ast\ast}(t)\right)  ^{\frac{j}{p}}\right)  .
\end{align*}
Therefore, we see that
\begin{align*}
t\left(  f_{(p)}^{\ast\ast}(t)-f_{(p)}^{\ast}(t)\right)   &  \leq\left\|
\left[  f-f^{\ast}(t)\right]  _{+}\right\|  _{p}^{p}+t\left(  \sum_{j=1}%
^{k}\left(
\begin{array}
[c]{c}%
p\\
j
\end{array}
\right)  f^{\ast}(t)^{p-j}\phi(t)^{j}\left(  \left|  \nabla f\right|
_{(p)}^{\ast\ast}(t)\right)  ^{\frac{j}{p}}\right) \\
&  \leq t\phi(t)^{p}\left|  \nabla f\right|  _{(p)}^{\ast\ast}(t)+t\left(
\sum_{j=1}^{k}\left(
\begin{array}
[c]{c}%
p\\
j
\end{array}
\right)  f^{\ast}(t)^{p-j}\phi(t)^{j}\left(  \left|  \nabla f\right|
_{(p)}^{\ast\ast}(t)\right)  ^{\frac{j}{p}}\right)  \text{\ \ \ \ (by
(\ref{Sk}))}\\
&  =t\phi_{(p)}(t)\left(  \left|  \nabla f\right|  _{(p)}^{\ast\ast}%
(t)+\sum_{j=1}^{k}\left(
\begin{array}
[c]{c}%
p\\
j
\end{array}
\right)  \left(  \frac{f_{(p)}^{\ast}(t)}{\phi_{(p)}(t)}\right)  ^{\frac
{p-j}{p}}\left(  \left|  \nabla f\right|  _{(p)}^{\ast\ast}(t)\right)
^{\frac{j}{p}}\right)  .
\end{align*}
Consequently,
\begin{equation}
\frac{f_{(p)}^{\ast\ast}(t)-f_{(p)}^{\ast}(t)}{\phi_{(p)}(t)}\leq\left|
\nabla f\right|  _{(p)}^{\ast\ast}(t)+\sum_{j=1}^{k}\left(
\begin{array}
[c]{c}%
p\\
j
\end{array}
\right)  \left(  \frac{f_{(p)}^{\ast}(t)}{\phi_{(p)}(t)}\right)  ^{\frac
{p-j}{p}}\left(  \left|  \nabla f\right|  _{(p)}^{\ast\ast}(t)\right)
^{\frac{j}{p}}. \label{norma1}%
\end{equation}
We can rewrite (\ref{norma1}) as
\begin{align*}
\frac{f_{(p)}^{\ast\ast}(t)-f_{(p)}^{\ast}(t)}{\phi_{(p)}(t)}  &  \leq\left|
\nabla f\right|  _{(p)}^{\ast\ast}(t)+\sum_{j=1}^{k}\left(
\begin{array}
[c]{c}%
p\\
j
\end{array}
\right)  \left(  \frac{f_{(p)}^{\ast}(t)}{\phi_{(p)}(t)}\right)  ^{\frac
{p-j}{p}}\left(  \left|  \nabla f\right|  _{(p)}^{\ast\ast}(t)\right)
^{\frac{j}{p}}+\frac{f_{(p)}^{\ast}(t)}{\phi_{(p)}(t)}-\frac{f_{(p)}^{\ast
}(t)}{\phi_{(p)}(t)}\\
&  =\left(  2^{\frac{k+1}{p}-1}\left(  \left|  \nabla f\right|  _{(p)}%
^{\ast\ast}(t)\right)  ^{1/p}+\left(  \frac{f_{(p)}^{\ast}(t)}{\phi_{(p)}%
(t)}\right)  ^{1/p}\right)  ^{p}-\frac{f_{(p)}^{\ast}(t)}{\phi_{(p)}(t)}\text{
\ \ (by (\ref{des2}))}%
\end{align*}
Hence
\[
\frac{f_{\left(  p\right)  }^{\ast\ast}(t)}{\phi_{(p)}(t)}\leq\left(
2^{\frac{k+1}{p}-1}\left(  \left|  \nabla f\right|  _{(p)}^{\ast\ast
}(t)\right)  ^{1/p}+\left(  \frac{f_{(p)}^{\ast}(t)}{\phi_{(p)}(t)}\right)
^{1/p}\right)  ^{p},
\]
yielding
\[
\left(  \frac{f_{\left(  p\right)  }^{\ast\ast}(t)}{\phi_{(p)}(t)}\right)
^{1/p}\leq2^{\frac{k+1}{p}-1}\left(  \left|  \nabla f\right|  _{(p)}^{\ast
\ast}(t)\right)  ^{1/p}+\left(  \frac{f_{(p)}^{\ast}(t)}{\phi_{(p)}%
(t)}\right)  ^{1/p}.
\]
Summarizing, we have obtained
\[
\left(  \frac{f_{\left(  p\right)  }^{\ast\ast}(t)}{\phi_{(p)}(t)}\right)
^{1/p}-\left(  \frac{f_{(p)}^{\ast}(t)}{\phi_{(p)}(t)}\right)  ^{1/p}%
\leq2^{\frac{k+1}{p}-1}\left(  \left|  \nabla f\right|  _{(p)}^{\ast\ast
}(t)\right)  ^{1/p}.
\]
$2\rightarrow3.$ Once again we use the elementary inequality
\[
\left(  x^{p}-y^{p}\right)  \leq p\left(  x-y\right)  \left(  x^{p-1}%
+y^{p-1}\right)  ,\text{ \ \ \ (}x\geq y\geq0),
\]
with $x=\left(  f_{\left(  p\right)  }^{\ast\ast}(t)\right)  ^{1/p}$ and
$y=\left(  f_{(p)}^{\ast}(t)\right)  ^{1/p}$. We obtain,
\begin{align*}
f_{\left(  p\right)  }^{\ast\ast}(t)-f_{(p)}^{\ast}(t)  &  \leq p\left(
\left(  f_{\left(  p\right)  }^{\ast\ast}(t)\right)  ^{1/p}-\left(
f_{(p)}^{\ast}(t)\right)  ^{1/p}\right)  \left(  \left(  f_{(p)}^{\ast\ast
}(t)\right)  ^{\frac{p-1}{p}}+\left(  f_{(p)}^{\ast}(t)\right)  ^{\frac
{p-1}{p}}\right) \\
&  \leq p2^{\frac{k+1}{p}-1}\phi(t)\left(  \left|  \nabla f\right|
_{(p)}^{\ast\ast}(t)\right)  ^{1/p}\left(  \left(  f_{(p)}^{\ast\ast
}(t)\right)  ^{\frac{p-1}{p}}+\left(  f_{(p)}^{\ast}(t)\right)  ^{\frac
{p-1}{p}}\right)  \text{ \ by (\ref{norma})}\\
&  \leq p2^{\frac{k+1}{p}-1}\phi(t)\left(  \left|  \nabla f\right|
_{(p)}^{\ast\ast}(t)\right)  ^{1/p}\left(  2\left(  f_{(p)}^{\ast\ast
}(t)\right)  ^{\frac{p-1}{p}}\right)  .
\end{align*}

Consequently,
\[
\frac{1}{p}\left(  f_{(p)}^{\ast\ast}(t)\right)  ^{\frac{1}{p}-1}\left(
f_{\left(  p\right)  }^{\ast\ast}(t)-f_{(p)}^{\ast}(t)\right)  \leq
p2^{\frac{k+1}{p}}\phi(t)\left(  \left|  \nabla f\right|  _{(p)}^{\ast\ast
}(t)\right)  ^{\frac{1}{p}}.
\]
Now observe that
\[
\frac{1}{p}\left(  f_{(p)}^{\ast\ast}(t)\right)  ^{\frac{1}{p}-1}\left(
\frac{f_{\left(  p\right)  }^{\ast\ast}(t)-f_{(p)}^{\ast}(t)}{t}\right)
=-\frac{\partial}{\partial t}\left(  \frac{1}{t}\int_{0}^{t}\left(  f^{\ast
}(s)\right)  ^{p}ds\right)  ^{1/p}.
\]
$3\rightarrow1.$

Let $\Omega\subset\subset X,$ and let $f\in Lip_{0}(\Omega),$ then, for
$t=\mu(\Omega),$ we have
\[
f_{\left(  p\right)  }^{\ast\ast}(t)=\frac{1}{t}\int_{0}^{t}\left(  f^{\ast
}(t)\right)  ^{p}dt=\frac{1}{t}\left\|  f\right\|  _{p}^{p}%
\]
and, similarly,
\[
\left|  \nabla f\right|  _{(p)}^{\ast\ast}(t)=\frac{1}{t}\left\|  \left|
\nabla f\right|  \right\|  _{p}^{p}.
\]
Since
\[
f_{(p)}^{\ast}(\mu(\Omega))=\inf_{x\in\Omega}\left|  f(x)\right|  ^{p}=0,
\]
the inequality (\ref{norma01}) becomes
\[
\frac{1}{t}\left\|  f\right\|  _{p}^{p}\leq p2^{\frac{k+1}{p}}\phi(\mu
(\Omega))\frac{1}{t}\left\|  \left|  \nabla f\right|  \right\|  _{p}\left\|
f\right\|  _{p}^{p-1},
\]
which is (\ref{desigualdad01}), up to constants.
\end{proof}

To complete the proof it remains to prove Lemma \ref{des}.

\begin{proof}
(of Lemma \ref{des}) We prove (\ref{des1}). Towards this end let us define
\[
f(x)=(x-b)^{p}-x^{p}+b^{p}+\sum_{j=1}^{k}\left(
\begin{array}
[c]{c}%
p\\
j
\end{array}
\right)  b^{p-j}(x-b)^{j},\text{ \ \ \ \ (}x\geq b).
\]
An elementary computation shows that $f(b)=\frac{\partial}{\partial
x}f(b)=\frac{\partial^{k-1}}{\partial x}f(b)=0.$ Moreover, since
\[
\frac{\partial^{k}}{\partial x}f(x)=p(p-1)\ldots(p-k+1)\left(  (x-b)^{p-k}%
-x^{p-k}+b^{p-k}\right)  ,
\]
and $0<p-k\leq1,$ we see that
\[
\left(  x-b\right)  ^{p-k}-x^{p-k}+b^{p-k}\geq0,
\]
consequently,
\[
f(x)\geq f(b)=0.
\]
To see (\ref{des2}) let us write $a=xb$ $(x\geq1).\ $We would like to show
that
\[
g(x)=(c(p)x+1)^{p}-x^{p}-1-\sum_{j=1}^{k}\left(
\begin{array}
[c]{c}%
p\\
j
\end{array}
\right)  (x-1)^{j}\geq0.
\]
An easy computation shows that $g(1)\geq0,\frac{\partial}{\partial x}%
g(1)\geq0,\cdots\frac{\partial^{k-1}}{\partial x}g(1)\geq0,$ and
$\frac{\partial^{k}}{\partial x}g(1)\geq0.$ Therefore, it will be enough to
prove that $\frac{\partial^{k+1}}{\partial x}g(x)\geq0.$ Again, by
computation, we find that
\[
\frac{\partial^{k+1}}{\partial x}g(x)=p(p-1)\ldots(p-k+1)(p-k)\left(
c(p)^{k+1}\left(  c(p)x+1\right)  ^{p-k-1}-x^{p-k-1}\right)  .
\]
Therefore the desired result will follow if we show that
\[
c(p)^{k+1}\left(  c(p)x+1\right)  ^{p-k-1}-x^{p-k-1}\geq0.
\]
Since $p-k-1<0,$ this amounts to show
\[
\frac{c(p)^{k+1}}{\left(  c(p)x+1\right)  ^{k+1-p}}\geq\frac{1}{x^{k+1-p}%
}\Leftrightarrow\frac{c(p)^{\frac{k+1}{k+1-p}}}{c(p)x+1}\geq\frac{1}%
{x}\Leftrightarrow xc(p)\left(  c(p)^{\frac{k+1}{k+1-p}}-1\right)  \geq1.
\]
But since
\[
c(p)^{\frac{k+1}{k+1-p}}-1\geq1\Leftrightarrow c(p)\geq2^{\frac{k+1}{p}-1},
\]
the desired result follows.
\end{proof}

\section{Final Remarks\label{secrem}}

In this section we show the explicit connection of our rearrangement
inequalities with the classical Sobolev inequalities, the Nash and Faber-Krahn
inequalities and point out possible directions for future research. In
particular, using Coulhon inequalities we will show a direct approach to some
self-improving properties of Sobolev inequalities for $p>1$.

\subsection{Nash Inequalities}

We start by giving a rearrangement characterization of the Nash type
inequalities. It was shown in \cite{bakr} (cf. also \cite{cou1}), that the
$(S_{\phi}^{p})$ conditions are equivalent to Nash type inequalities. As a
consequence, the results of this paper give a characterization of Nash
inequalities in terms of rearrangements which we shall now describe.

We first observe that, with some trivial changes, one can adapt the proof of
Proposition 2.4 in \cite{cou1} (case $p=2)$ to obtain the following
equivalence (for Nash inequalities for $p>1)$

\begin{proposition}
Let $p>1.$ The following are inequalities are equivalent up to multiplicative constants

(i) $(S_{\phi}^{p})$ holds

(ii) There exist positive constants $c_{1}$ and $c_{2}$ such that
\[
\left\|  f\right\|  _{p}\leq c_{1}\phi\left(  c_{2}\left(  \frac{\left\|
f\right\|  _{1}}{\left\|  f\right\|  _{p}}\right)  ^{\frac{p}{p-1}}\right)
\left\|  \left|  \nabla f\right|  \right\|  _{p}%
\]
for all $f\in Lip_{0}(X)$.
\end{proposition}

The case $\phi(t)=t^{1/n}$, $p=2,$ corresponds to the classical Nash
inequality
\[
\left\|  f\right\|  _{2}^{1+2/n}\leq c\left\|  f\right\|  _{1}^{2/n}\left\|
\left|  \nabla f\right|  \right\|  _{2}.
\]
Therefore, by Theorem \ref{teo}, Nash's inequality is equivalent to
\[
\left(  \frac{f_{\left(  2\right)  }^{\ast\ast}(t)}{t^{2/n}}\right)
^{1/2}-\left(  \frac{f_{(2)}^{\ast}(t)}{t^{2/n}}\right)  ^{1/2}\preceq\left(
\left|  \nabla f\right|  _{(2)}^{\ast\ast}(t)\right)  ^{1/2},\text{ }f\in
Lip_{0}(\mathbb{R}^{n}).
\]

\subsection{Classical Sobolev Inequalities}

We now consider a new approach, via rearrangement inequalities, of the known
(cf. \cite{cou}, \cite{bakr}, \cite{cou1} and the references therein)
equivalence between the classical Euclidean Sobolev inequalities and Coulhon
inequalities. The case $p=1$ of (\ref{norma}) gives us the inequality
(\ref{berta}), whose connection to Sobolev inequalities was discussed
extensively elsewhere (cf. \cite{mamiadv}).

Let us consider the case $1\leq p<n,$ $\frac{1}{\bar{p}}=\frac{1}{p}-\frac
{1}{n}.$ Let $\phi(t)=t^{1/n}.$ We shall denote the corresponding $(S_{\phi
}^{p})$ condition by $(S_{n}^{p}).$ Our aim is to prove that $(S_{n}^{p})$
implies the classical Sobolev inequality
\[
\left\|  f\right\|  _{L(\bar{p},p)}\preceq\left\|  \nabla f\right\|  _{L^{p}%
},f\in Lip_{0}(\mathbb{R}^{n}),
\]
where for $1\leq r<\infty,1\leq q\leq\infty,$%
\[
\left\|  f\right\|  _{L(r,q)}=\left\{  \int_{0}^{\infty}\left(  f^{\ast
}(t)t^{\frac{1}{r}}\right)  ^{q}\frac{dt}{t}\right\}  ^{1/q}.
\]
By a well known result, apparently originally due to Maz'ya, weak type Sobolev
inequalities self-improve to strong type Sobolev inequalities (cf.
\cite{bakr}, \cite{mmp}, and the references therein). We shall discuss this
self-improvement in detail in the next subsection. Taking this fact for
granted, it will be enough to show that $(S_{n}^{p})$ implies the weak type
Sobolev inequality
\begin{equation}
\left\|  f\right\|  _{L(\bar{p},\infty)}\preceq\left\|  \nabla f\right\|
_{L^{p}},f\in Lip_{0}(\mathbb{R}^{n}), \label{tio}%
\end{equation}
where
\[
\left\|  f\right\|  _{L(\bar{p},\infty)}=\sup_{t}\{f^{\ast}(t)t^{1/\bar{p}%
}\}.
\]
To prove (\ref{tio}) let us first recall that since $\bar{p}>1,$ for $f\in
Lip_{0}(\mathbb{R}^{n})$ we have (cf. \cite{bs}, \cite{bmr}),
\[
\left\|  f\right\|  _{L(\bar{p},\infty)}\simeq\sup_{t}\{\left(  f^{\ast\ast
}(t)-f^{\ast}(t)\right)  t^{1/\bar{p}}\}.
\]
We have shown above that $(S_{n}^{p})$ implies (\ref{norma}); therefore it
follows that
\begin{align*}
\left(  f_{p}^{\ast\ast}(t)\right)  ^{1/p}-f^{\ast}(t)  &  \preceq
t^{1/n}\left(  \left|  \nabla f\right|  _{p}^{\ast\ast}(t)\right)  ^{1/p}\\
&  =t^{1/n-1/p}\left\{  \int_{0}^{t}\left|  \nabla f\right|  ^{\ast}%
(s)^{p}ds\right\}  ^{1/p}.
\end{align*}
Combining the last inequality with Jensen's inequality we get
\begin{align}
f^{\ast\ast}(t)-f^{\ast}(t)  &  \leq\left(  f_{p}^{\ast\ast}(t)\right)
^{1/p}-f^{\ast}(t)\nonumber\\
&  \preceq t^{1/n-1/p}\left\{  \int_{0}^{t}\left|  \nabla f\right|  ^{\ast
}(s)^{p}ds\right\}  ^{1/p}. \label{abaco}%
\end{align}
Summarizing, for $f\in Lip_{0}(\mathbb{R}^{n}),$%
\begin{align*}
\left\|  f\right\|  _{L(\bar{p},\infty)}  &  \simeq\sup_{t}\{\left(
f^{\ast\ast}(t)-f^{\ast}(t)\right)  t^{1/\bar{p}}\}\\
&  \preceq\sup_{t}\left\{  \int_{0}^{t}\left|  \nabla f\right|  ^{\ast}%
(s)^{p}ds\right\}  ^{1/p}\\
&  \leq\left\|  \left|  \nabla f\right|  \right\|  _{p},
\end{align*}
as we wished to show.

In this next section we shall discuss in detail the case $p=n,$ and show the
self-improvement of Sobolev-Coulhon inequalities.

\subsection{Self-improvement}

There are several known mechanisms to show the self-improvement of Sobolev
inequalities. Here we choose to adapt a variant the method apparently first
developed by Maz'ya-Talenti (cf. \cite{mamiadv} for a generalized version)
using differential inequalities, focussing on the Euclidean case.

For a domain $\Omega\subset\mathbb{R}^{n},$ we have (cf. \cite{mmp} for the
classical Euclidean case or \cite{mamiadv} for the general metric space case)
the following formulation of the Polya-Szeg\"{o} principle
\begin{equation}
\left(  \int_{0}^{\left|  \Omega\right|  }\left(  s^{1-\frac{1}{n}}\left(
-f^{\ast}\right)  ^{^{\prime}}(s)\right)  ^{p}ds\right)  ^{1/p}\preceq\left(
\int_{0}^{\left|  \Omega\right|  }\left(  \left|  \nabla f\right|  ^{\ast
}(s)\right)  ^{p}ds\right)  ^{1/p},\text{ }p\geq1,\text{ }f\in Lip_{0}%
(\Omega). \label{negada}%
\end{equation}

To use this powerful inequality we now reformulate (\ref{abaco}) as an
elementary differential inequality. For $f\in Lip_{0}(\Omega),$ let
$F(t):=\left(  f^{\ast\ast}(t)-f^{\ast}(t)\right)  ^{p}$ $t^{1-\frac{p}{n}%
},1\leq p<n.$ Then $F$ is a positive, absolutely continuous function (cf.
\cite{leoni}), which by (\ref{abaco}) satisfies
\[
F(t)\preceq\int_{0}^{t}\left(  \left|  \nabla f\right|  ^{\ast}(s)\right)
^{p}ds.
\]
It follows that $F(0)=0,$ and therefore we can write $F(t)=\int_{0}%
^{t}F^{^{\prime}}(s)ds,t>0.$ We estimate $F$ through this representation. By
direct computation,
\begin{align*}
F^{^{\prime}}(t)  &  =(1-\frac{p}{n})t^{-\frac{p}{n}}[f^{\ast\ast}(t)-f^{\ast
}(t)]^{p}+t^{1-\frac{p}{n}}p[f^{\ast\ast}(t)-f^{\ast}(t)]^{p-1}\left[  \left(
f^{\ast\ast}(t)\right)  ^{\prime}-\left(  f^{\ast}\right)  ^{^{\prime}%
}(t)\right] \\
&  =(1-\frac{p}{n})t^{-\frac{p}{n}}[f^{\ast\ast}(t)-f^{\ast}(t)]^{p}%
+t^{1-\frac{p}{n}}p[f^{\ast\ast}(t)-f^{\ast}(t)]^{p-1}\left[  (-1)\left(
\frac{f^{\ast\ast}(t)-f^{\ast}(t)}{t}\right)  -\left(  f^{\ast}\right)
^{^{\prime}}(t)\right] \\
&  =(1-\frac{p}{n}-p)t^{-\frac{p}{n}}[f^{\ast\ast}(t)-f^{\ast}(t)]^{p}%
+t^{1-\frac{p}{n}}p[f^{\ast\ast}(t)-f^{\ast}(t)]^{p-1}\left(  -f^{\ast
}\right)  ^{^{\prime}}(t).
\end{align*}
The previous computation, combined with the fact that $F(t)$ is positive,
yields%
\[
(-1)(1-\frac{p}{n}-p)\int_{0}^{\left|  \Omega\right|  }[f^{\ast\ast
}(t)-f^{\ast}(t)]^{p}t^{-\frac{p}{n}}dt\leq p\int_{0}^{\left|  \Omega\right|
}t^{1-\frac{p}{n}}[f^{\ast\ast}(t)-f^{\ast}(t)]^{p-1}\left(  -f^{\ast}\right)
^{^{\prime}}(t)dt.
\]
H\"{o}lder's inequality and (\ref{negada}) yields%
\begin{align*}
&  \int_{0}^{\left|  \Omega\right|  }\left(  [f^{\ast\ast}(t)-f^{\ast
}(t)]t^{\frac{1}{\bar{p}}}\right)  ^{p}\frac{dt}{t}\\
&  =\int_{0}^{\left|  \Omega\right|  }t^{-\frac{p}{n}}[f^{\ast\ast}%
(t)-f^{\ast}(t)]^{p}dt\\
&  \leq\frac{p}{(\frac{p}{n}+p-1)}\int_{0}^{\left|  \Omega\right|  }%
t^{1-\frac{p}{n}}[f^{\ast\ast}(t)-f^{\ast}(t)]^{p-1}\left(  -f^{\ast}\right)
^{^{\prime}}(t)dt\\
&  =\frac{p}{(p-1+\frac{p}{n})}\int_{0}^{\left|  \Omega\right|  }([f^{\ast
\ast}(t)-f^{\ast}(t)]^{p-1}t^{\frac{1-p}{n}})(t^{1-\frac{1}{n}}\left(
f^{\ast}\right)  ^{^{\prime}}(t))dt\\
&  \preceq\left(  \int_{0}^{\left|  \Omega\right|  }\left(  t^{\frac{1-p}{n}%
}[f^{\ast\ast}(t)-f^{\ast}(t)]^{p-1}\right)  ^{\frac{p}{p-1}}dt\right)
^{(p-1)/p}\left(  \int_{0}^{\left|  \Omega\right|  }\left(  t^{1-\frac{1}{n}%
}\left(  -f^{\ast}\right)  ^{^{\prime}}(t)\right)  ^{p}dt\right)  ^{1/p}\\
&  =c_{n,p}\left(  \int_{0}^{\left|  \Omega\right|  }\left(  t^{\frac{1}%
{\bar{p}}}[f^{\ast\ast}(t)-f^{\ast}(t)]\right)  ^{p}\frac{dt}{t}\right)
^{1/p^{\prime}}\left\|  \left|  \nabla f\right|  \right\|  _{p}.
\end{align*}

Consequently, assuming \textit{apriori} that $\int_{0}^{\left|  \Omega\right|
}\left(  t^{\frac{1}{p}-\frac{1}{n}}[f^{\ast\ast}(t)-f^{\ast}(t)]^{p}\right)
\frac{dt}{t}<\infty,$ we see that
\[
\left\{  \int_{0}^{\left|  \Omega\right|  }\left(  [f^{\ast\ast}(t)-f^{\ast
}(t)]t^{\frac{1}{\bar{p}}}\right)  ^{p}\frac{dt}{t}\right\}  ^{1/p}%
\preceq\left\|  \left|  \nabla f\right|  \right\|  _{p}.
\]
For $f\in Lip_{0}(\Omega)$ all the formal calculations above can be easily
justified and we find the sharp Sobolev inequality
\[
\left\|  f\right\|  _{L(\bar{p},p)}\simeq\left\{  \int_{0}^{\left|
\Omega\right|  }\left(  [f^{\ast\ast}(t)-f^{\ast}(t)]t^{\frac{1}{\bar{p}}%
}\right)  ^{p}\frac{dt}{t}\right\}  ^{1/p}\preceq\left\|  \left|  \nabla
f\right|  \right\|  _{p}.
\]

Let us note that the previous calculation also works for $p=n.$ In this case
we should let $\frac{1}{\bar{p}}=0$ and we obtain%
\[
\left\{  \int_{0}^{\left|  \Omega\right|  }\left(  f^{\ast\ast}(t)-f^{\ast
}(t)\right)  ^{n}\frac{dt}{t}\right\}  ^{1/n}\preceq\left\|  \left|  \nabla
f\right|  \right\|  _{n}.
\]
In this case the left hand side should be re-interpreted as the *norm* of
$L(\infty,n)$, the space defined by the condition (cf. \cite{bmr})%
\[
\left\{  \int_{0}^{\left|  \Omega\right|  }\left(  f^{\ast\ast}(t)-f^{\ast
}(t)\right)  ^{n}\frac{dt}{t}\right\}  ^{1/n}<\infty.
\]
It was shown in \cite{bmr} that this condition implies the classical
exponential integrability results of Trudinger and Brezis-Wainger.

Note that the self-improvement for general $\phi,$ which we have not discussed
here, will involve the $p-$Lorentz $\Lambda_{\phi}$ spaces (for further
related discussions we refer to \cite{mamicon}).

\subsection{The Morrey-Sobolev theorem}

The connection between rearrangement inequalities and the Morrey-Sobolev
theorem (i.e. the case $p>n$ of the Sobolev embedding theorem$)$ has been
treated at great length in our recent article \cite{mamiarxiv}. We consider
here the corresponding Coulhon variant, but, once again for the sake of
brevity, and to avoid technical complications, we shall only sketch the
details for Sobolev spaces $W_{0}^{1}(Q)$ on the cube $Q=(0,1)^{n}$.

In this section we let $p>n,$ then $\frac{1}{\bar{p}}=\frac{1}{p}-\frac{1}%
{n}<0.$ Using the fact that $(-f^{\ast\ast}(t))^{\prime}=\frac{f^{\ast\ast
}(t)-f^{\ast}(t)}{t}$ we can integrate the inequality (\ref{abaco}) to obtain
\begin{align*}
f^{\ast\ast}(0)-f^{\ast\ast}(1)  &  =\int_{0}^{1}\frac{f^{\ast\ast}%
(t)-f^{\ast}(t)}{t}dt\\
&  \preceq\int_{0}^{1}t^{-\frac{1}{\bar{p}}}\left(  \int_{0}^{t}\left|  \nabla
f\right|  ^{\ast}(s)^{p}ds\right)  ^{1/p}\frac{dt}{t}\\
&  \leq\left\|  \left|  \nabla f\right|  \right\|  _{p}\int_{0}^{1}%
t^{-\frac{1}{\bar{p}}-1}dt\\
&  =c_{p}\left\|  \left|  \nabla f\right|  \right\|  _{p}.
\end{align*}
Extending the inequalities we have obtained in this note through the use of
signed rearrangements, and using an extension of a scaling argument that
apparently goes back to \cite{garsia} (we must refer to \cite[pag.
3]{mamiarxiv} for more details) we find that given $x,y\in Q,$
\begin{align*}
\left|  f(x)-f(y)\right|   &  \preceq\left\|  \left|  \nabla f\right|
\right\|  _{p}\left|  x-y\right|  ^{n(\frac{1}{n}-\frac{1}{\bar{p}})}\\
&  =\left\|  \left|  \nabla f\right|  \right\|  _{p}\left|  x-y\right|
^{1-\frac{n}{p}}.
\end{align*}

\subsection{Further connections}

In this section we mention some problems and possible projects we find of some interest.

In the literature there are other definitions of the notion of gradient in the
metric setting (e.g. \cite{haj} and the references therein) and it remains an
open problem to fully explore the connections with our development
here\footnote{For partial results (restricted to doubling measures) connecting
different notions of the gradient with rearrangement inequalities we refer to
\cite{badr}, \cite{kalis} and the references therein.}.

We hope to discuss the connection between isoperimetry, rearrangements and
discrete Sobolev inequalities elsewhere.

For aficionados of interpolation theory we should note that, while there are
obvious connections between the$\ (S_{\phi}^{p})$ conditions and the
$J-$method of interpolation or perhaps, even more appropriately, with the
corresponding version of this method for the $E-$method of approximation (cf.
\cite{jm}), we could not find a treatment in the literature. Such
considerations are somehow implicit in the approach given in \cite{bakr}, and
more explicitly in the unpublished manuscript \cite{cm}. Likewise, the $\phi$
inequalities that appear in the formulation of Nash's inequality
above\footnote{There is an extensive literature on $\phi$ inequalities (cf.
\cite{bcr}, and the references therein).}, appear directly related to the
$K/J$ inequalities of the extrapolation theory of \cite{jm}.

Still another direction for future research is to develop in more detail the
connection of the results in this paper and the work of Xiao \cite{xiao} on
the $p-$Faber-Krahn inequality.

Finally in this section we have discussed only the Euclidean case. It will be
of interest to develop a detailed treatment of these applications in the
general metric case.

\begin{acknowledgement}
We are grateful to the referee for many suggestions to improve the quality of
the paper. The research for this paper was partly carried out while the second
named author was visiting the Universidad Autonoma de Barcelona. He is
grateful to this institution and the CRM for excellent working conditions.
\end{acknowledgement}

\end{document}